\documentclass[sn-mathphys,Numbered]{sn-jnl}% Math and Physical Sciences Reference Style
\usepackage{graphicx}%
\usepackage{multirow}%
\usepackage{amsmath,amssymb,amsfonts}%
\usepackage{amsthm}%
\usepackage{mathrsfs}%
\usepackage[title]{appendix}%
\usepackage{xcolor}%
\usepackage{textcomp}%
\usepackage{manyfoot}%
\usepackage{booktabs}%
\usepackage{algorithm}%
\usepackage{algorithmicx}%
\usepackage{algpseudocode}%
\usepackage{listings}%
\newtheorem{rmk}{Remark}[section]
\newtheorem{prop}{Proposition}[section]
\newtheorem{lemma}{Lemma}[section]
\newtheorem{defn}{Definition}[section]
\newtheorem{exm}{Example}[section]
\newtheorem{quest}{Question}[section]
\newtheorem{corollary}{Corollary}[section]

\theoremstyle{thmstylethree}%

\begin{document}

\title[Article Title]{On Nil-semicommutative Modules}


\author[1]{\fnm{M.} \sur{Rhoades}}\email{rhoades.phd.math@manipuruniv.ac.in}

\author*[1]{\fnm{Kh.} \sur{Herachandra}}\email{heramath@manipuruniv.ac.in}

\author[2]{\fnm{Nazeer} \sur{Ansari}}\email{nazee.ansari@rgu.ac.in}

\affil[1]{\orgdiv{Department of Mathematics}, \orgname{Manipur University}, \city{Canchipur}, \postcode{795003}, \state{Manipur}, \country{India}}

\affil[2]{\orgdiv{Department of Mathematics}, \orgname{Rajiv Gandhi University}, \city{Papum Pare}, \postcode{791112}, \state{Arunachal Pradesh}, \country{India}}

%\affil[3]{\orgdiv{Department}, \orgname{Organization}, \orgaddress{\street{Street}, \city{City}, \postcode{610101}, \state{State}, \country{Country}}}


\abstract{	In this paper, we introduce a new concept in Nil-semicommutative modules and present it as an extension of  Nil-semicommutative rings to modules. We prove that the class of Nil-semicommutative modules is contained in the class of Weakly semicommutative modules while that of the converse may not be true. We also show that in case of Semicommutative modules and Nil-semicommutative modules, one does not imply the other. Moreover, for a given Nil-semicommutative ring, we provide the conditions under which the same can be extended to a Nil-semicommutative module. Lastly, we also prove that for a left $R$-module $M$, $_RM$ is Nil-semicommutative iff it's localization $S^{-1}M$ over the ring $S^{-1}R$ is also Nil-semicommutative.Various other examples and propositions highlighting the comparative studies of this new class of modules with different classes of modules are also discussed in order to validate the concept.}




\maketitle

\section{Introduction}\label{sec1}
	Throughout the study, the ring $R$ is associative with identity, and the module $M$ is unitary left $R$-module.~$M_{n}(R)$, $T_n(R)$ and $S_n(R)$ denote respectively the ring of $n \times n$ matrices, upper triangular matrices and special upper triangular matrices with the same diagonal entries over $R$. An element $m$ is torsion if ~$\exists~ 0 \neq r \in R$ such that $rm = 0$ and the collection of such elements is denoted by   $Tor(M)$.~A torsion-free module $M$ is a module satisfying the condition $Tor(M) = 0$. For $S_0$ denoting the set of non-zero zero divisors of a ring $R$, $T(M) = \{ m \in M : rm = 0$ for some $0 \neq r \in S_0 \} $. $Nil(R)$ denotes the set of all nilpotent elements in the ring $R$ and $Nil_R (M)$ denotes the set of all nilpotent elements in the left $R$-module $M$.\\

The concept of $Nil(R)$ has been there for long, and the utilization of this concept by many algebraists can be seen in numerous concepts in the form of  Reduced rings, Semicommutative rings, Weakly Semicommutative rings, Nil-semicommutative rings, etc.~In fact, it has been studied and proved in \textcolor{blue}{\citep{chen}} that Semicommutative ring $\Rightarrow$ Nil-semicommutative ring $\Rightarrow$ Weakly Semicommutative ring with the converses not holding in each of the cases. Several other results based on the concept of nilpotency in rings are also being studied and discussed presently.\\

However, the concept of nilpotent elements in a module was defined and studied by Ssevviiri and Groenewald \textcolor{blue}{\cite{ssev}} in 2013 and with that a breakthrough in extending existing results on nilpotency in rings to their modules counterparts was achieved.~An element $m \in M$ is said to be a nilpotent element if either $m = 0$ or if for $m \neq 0$, $\exists$ $r \in R$ and $k \in \mathbb{N}$ such that $r^km = 0$ but $rm \neq 0$.\\
%E.g.~i) ~~$\bar0$,$\bar1$,$\bar3$ are nilpotent elements in the $Z$-module $Z_{4}$ \\
%E.g. ii) ~~$\bar1$ is nilpotent in the $Z$-module $Z_{n^k}$ for $k > 1$\\

Moreover, using this definition, Ansari and Herachandra in \textcolor{blue}{\cite{1}} came up with a result that further simplified the process of finding nilpotent elements in modules. In their paper, they proved the following proposition.\\

Let $m$ be an element of the left $R$-module $M$. Then, the following conditions are equivalent.
\begin{itemize}
	\item[(i)]  $\exists~ r \in R$ and $n \geq 2$ such that $r^{n}m = 0 ; r^{n-1}m \neq 0$. 
	\item[(ii)]  $\exists~ t \in R$ such that $t^{2}m = 0 ; tm \neq 0$. 
\end{itemize}

So, if one wants to check whether an element $m \in Nil_R(M)$, it suffices to check if either one of the above conditions holds or not.\\ %However, it is noted that the converse of this does not hold as $ \exists~ m \in Nil_R(M) $, which does not satisfy the condition (b) above (take $m = 0$).\\

Revisiting the journey so far, it is learnt that after an inspiration from a result on Reduced rings~(rings with no non-zero nilpotent elements) in \cite{armn} by E. Armendariz, Rege, and Chhawchharia introduced the concept of Armendariz rings in \textcolor{blue}{\cite{rc}}. A ring $R$ is Armendariz if given any polynomials $f(x) =\sum_{i=0}^{n}a_ix^i $ and $g(x) = \sum_{j=0}^{m}b_jx^j$ with coefficients in $R$ satisfying $f(x).g(x) = 0$, we have $a_ib_j = 0 ~\forall~ i, j$.~Further, many new studies were made on this for some time, and in \textcolor{blue}{\cite{ant}}, Antoine introduced the concept of Nil Armendariz rings as a generalization of Armendariz rings. It is defined that a ring $R$ is Nil Armendariz if for polynomials $f(x) = \sum_{i=0}^{n}a_ix^i $ and $g(x) = \sum_{j=0}^{m}b_jx^j$ with co-efficients in $R$ satisfying the condition $f(x).g(x) \in Nil(R)[x]$, we have $a_ib_j \in Nil(R)~ \forall i, j$. Then, in \textcolor{blue}{\cite{liu}}, Liu and Zhao introduced Weak Armendariz rings and proved it as a further generalization of Nil Armendariz rings.~In their words, a ring $R$ is Weak Armendariz ring if for polynomials $f(x) =\sum_{i=0}^{n}a_ix^i $ and $g(x) =\sum_{j=0}^{m}b_jx^j $ with coefficients in $R$ satisfying the condition $f(x).g(x) = 0$, we have $a_ib_j \in Nil(R)$. Thus, in a nutshell, one can easily see that Armendariz rings $\Rightarrow$ Nil Armendariz rings $\Rightarrow$ Weak Armendariz rings.   \\

Similarly, in \textcolor{blue}{\cite{lee}}, Lee and Zhou introduced the concept of Reduced modules and also provided numerous generalizations in the module counterparts of the previously mentioned rings. In \textcolor{blue}{\cite{lee}}, a left $R$-module $M$ is reduced if it satisfies any of the following conditions:

\begin{itemize}
	\item[(i)]
	whenever $a \in R, m \in M$ satisfy $a^2m = 0$, we have $aRm = 0$
	\item[(ii)] whenever $a \in R, m \in M$ satisfy $am = 0$, we have $aM \bigcap Rm = 0$
\end{itemize}

%(a) whenever $a \in R, m \in M$ satisfy $a^2m = 0$, we have $aRm = 0$\\
%(b) whenever $a \in R, m \in M$ satisfy $am = 0$, we have $aM \ Rm = 0$\\

Moreover, the concept of Armendariz modules was also introduced in \textcolor{blue}{\cite{rc}} and numerous studies on the relationship of this class of modules with the class of Reduced modules and Semicommutative modules were done by Rege and Buhphang in \textcolor{blue}{\cite{br}}.~In \textcolor{blue}{\cite{rc}} a left $R$-module $M$ is said to be Armendariz if whenever $f(x) = \sum_{i=0}^{n}a_ix^i \in R[x]$ and $m(x) =\sum_{j=0}^{m}m_jx^j  \in M[x]$ satisfy $f(x).m(x) = 0$, then $a_im_j = 0 ~\forall~ i,j$.\\

Then, Ansari and Herachandra in \textcolor{blue}{\cite{2}} made further studies on the concept of nilpotent elements in modules and their role and impact on various types of modules. In fact, it was proved in \textcolor{blue}{\cite{2}} that for a reduced module, there are absolutely no non-zero nilpotent elements. Further, the extension of the concept of Weak Armendariz rings to Weak Armendariz modules was also done in \textcolor{blue}{\cite{2}}, thereby generalizing the said ring theoretical concept. A left $R$-module $M$ is said to be Weak Armendariz if whenever $f(x) = \sum_{i=0}^{n}a_ix^j\in R[x]$ and $m(x) = \sum_{j=0}^{k}m_jx^j\in M[x]$ satisfy $f(x)m(x) = 0$, we have $a_im_j \in Nil_R(M) ~\forall~ i,j$.\\   

On the other hand, Semicommutative rings were studied in \textcolor{blue}{\cite{hls}} by Huh et. al in 2002, and various other studies on it were also done under the names of IFP rings and zero-insertive(ZI) rings in the literature. A generalization on the same, Semicommutative modules was then introduced and studied by Buhphang and Rege \textcolor{blue}{\cite{br}} in the same year.~In \textcolor{blue}{\cite{1}} , it is defined that a ring is said to be Semicommutative if whenever
$ a, b \in R$ satisfy $ab = 0$, we have $aRb = 0$.~And in in \textcolor{blue}{\cite{br}}, it is defined that a left $R$-module $M$ is said to be Semicommutative module if whenever 
$a \in  R$ and $m \in M$ satisfy $am = 0$, we have $aRm = 0.$\\

Weakly semicommutative rings were introduced in \textcolor{blue}{\cite{liang}} and numerous studies were made on it by Liang et.al. A ring is said to be Weakly semicommutative if whenever 
$a,b \in  R$ satisfy $ab = 0$, we have $aRb \subseteq  Nil(R).$ Inspired by this, in \textcolor{blue}{\cite{1}} Ansari and Herachandra introduced Weakly Semicommutative modules as a generalization of Weakly semicommutative rings. A left $R$ module $M$ is said to be Weakly semicommutative if whenever 
$a \in R$ and $m \in M$ satisfy $am = 0$, then $aRm \subseteq Nil_R (M)$.\\ 

Lastly, Nil-semicommutative rings were defined and studied by Chen \textcolor{blue}{\cite{chen}} in 2011, with further contributions made by Moussavi et. al in 2012. In \textcolor{blue}{\cite{chen}}, a ring is said to be Nil-semicommutative if whenever
$a, b \in  R$ satisfy $ab \in Nil(R)$, we have $aRb \subseteq Nil(R).$\\\\

\section{Nil-semicommutative module}\label{sec2}

In this section, we introduce Nil-semicommutative module and investigate its preliminary results. We also show that many properties of Semicommutative modules can be extended to that of Nil-semicommutative modules under the general setting. \\

\begin{defn} A left $R$- module $M$ is said to be Nil semicommutative if whenever 
	$a \in R$ and $m \in M$ satisfy $am \in Nil_R (M)$, we have $aRm \subseteq Nil_R (M).$\\
	
\end{defn}

Recall that in the case of modules, Reduced $\Rightarrow$ Semicommutative $\Rightarrow$ Weakly semicommutative. Similarly, by definition of Nil-semicommutative modules, it is clear to say that Reduced $\Rightarrow$ Nil-semicommutative $\Rightarrow$  Weakly semicommutative. However, as discussed in  Example \ref{eg2},~the converse may not be true in general. Moreover, in the case of rings, we know that every Semicommutative ring is Nil-semicommutative, and it looks like the same follows in the case of modules. But, the Example \ref{eg2} given below will eliminate such possibilities, i.e., there exist modules that are Nil-semicommutative but not Semicommutative, and there are modules that are Semicommutative but not Nil-semicommutative.\\

\begin{rmk} A Weakly semicommutative module may not necessarily be a Nil-semicommutative module. The following  Example \ref{eg2} is one in this regard.\\ 
\end{rmk}

%\begin{exm}\label{eg1}
%	Consider the $\mathbb{Z}$-module ${\mathbb{Z}}_{p^n}$, then, ${\mathbb{Z}}_{p^n}$ is a \it{Weakly semicommutative} module for $n\geq 2$ which is not Nil-semicommutative module. This is because ${\mathbb{Z}}_{p^n}$ being a module over commutative ring is always \it{Semicommutative} and hence, \it{Weakly semicommutative}.\\
%\end{exm}

\begin{rmk}\label{rmk1}
	A Semicommutative module may not necessarily be a Nil-semicommutative module. The following example  Example \ref{eg2} is one in this regard.\\
\end{rmk}

\begin{lemma}\label{lem1}
	For the left $\mathbb{Z}_n$- module $\mathbb{Z}_n$, $\bar{1} \in Nil_{\mathbb{Z}_n}(\mathbb{Z}_n)$ $\iff$ $n$ is not square free.\\
\end{lemma}

\noindent\textbf{Proof:}  Suppose, $n$ is square free. Then, $ n = p_1p_2p_3...p_k$ ~for primes $p_i$ and $k \in \mathbb{N}$. Now, as $ \bar{1} \in Nil_{\mathbb{Z}_n}(\mathbb{Z}_n) ~\exists~ \bar{r} \in \mathbb{Z}_n$ and $k_0 \in \mathbb{N}$ ~such that~ $ \bar{r}^{k_0}.\bar{1} = \bar{0}$ ~but~ $ \bar{r}.\bar{1} \neq \bar{0}$. This means that $ r^{k_0} = n.l$ for some $ l \in \mathbb{Z}$ and so, $ p_i ~|~ r^{k_0}  ~\forall~ i = 1$ to $k$.~Consequently,  $ p_i ~|~ r ~\forall~ i$ and we have $ r= (p_1p_2p_3...p_k).t$ ~for some~ $ t \in \mathbb{Z}$. However, this means that $ \bar{r}.\bar{1} = \bar{0}$ which is a contradiction. Hence, $n$ cannot be square free.

For the converse, since $n$ is not square free, $ n = p_1^k.s$ for $ k \geq 2$ and some $s \in \mathbb{Z}$. Now, $ \exists~ (\bar{p_1.s}) \in \mathbb{Z}_n$ such that $ (\bar{p_1.s})^k.\bar{1} = \bar{0}$ ~but~ $ (\bar{p_1.s}).\bar{1} \neq \bar{0}$. Hence, $ \bar{1} \in Nil_{\mathbb{Z}_n}(\mathbb{Z}_n) $.\\

\begin{exm}\label{eg2}
	Consider the $\mathbb{Z}_{p^n}$-module $\mathbb{Z}_{p^n}$. Then, $\mathbb{Z}_{p^n}$ being a module over a commutative ring is Semicommutative, and also Weakly Semicommutative. However, for $n \geq 2$, $\mathbb{Z}_{p^{n}}$ is not  Nil-semicommutative.\\
\end{exm}

\noindent\textbf{Proof:}  By Lemma \ref{lem1}, we know that $ \bar{1}.\bar{1} \in Nil_{\mathbb{Z}_{p^n}}(\mathbb{Z}_{p^n})$ [since $ p^n$ is not square free for $ n\geq 2$]. But, we claim $ \bar{1}.\bar{p}^{(n-1)}.\bar{1} \notin Nil_{\mathbb{Z}_{p^n}}(\mathbb{Z}_{p^n})$. Suppose this fails, i.e. $ \bar{p}^{(n-1)} \in Nil_{\mathbb{Z}_{p^n}}(\mathbb{Z}_{p^n})$. Then, $ ~\exists~ \bar{r} \in \mathbb{Z}_{p^n}$ and $ k_0 \in \mathbb{N}$ such that  $ \bar{r}^{k_0}.\bar{p}^{(n-1)} = \bar{0}$ ~but~ $ \bar{r}.\bar{p}^{(n-1)} \neq \bar{0}$. This means that $r^{k_0}.p^{(n-1)} = p^n.s$ for some $ s \in \mathbb{Z}$. Consequently, $ p ~|~ r^{k_0}$ and hence, $p~|~r$. So, we have $ r = p.l$ for some $ l \in \mathbb{Z}$. This means that $ \bar{r}\bar{p}^{n-1} = \bar{pl}.\bar{p}^{n-1} = \bar{0}$, which is a contradiction. Therefore, $ \bar{1}.\bar{p}^{(n-1)}.\bar{1} \notin Nil_{\mathbb{Z}_{p^n}}(\mathbb{Z}_{p^n})$. Hence, $\mathbb{Z}_{p^n}$ is not Nil-semicommutative.\\

\begin{lemma}\label{lem2}
	Let $M$ be a left $R$-module.\\ Then for $ n\geq 2, Nil_{M_n(R)}(M_{n}(M))= ~ _{M_n(R)}M_{n}(M)$.\\
\end{lemma}

\noindent	\textbf{Proof:}
Consider any  non-zero matrix $K = [m_{ij}]_{n\times n}$. This implies atleast one $m_{ij}\neq 0$ for some $1 \leq i,j\leq n$. Then, two cases arise:
\begin{itemize}
	\item[(i)] For $i\neq j$, we can take $r=e_{ji}$. Thus we can easily see that $r^{2}K=(e_{ji})^{2}K=0$, but $rK=e_{ji}K \neq 0$.
	\item[(ii)] For $i=j$, we can take $r=e_{li}$ such that $l \neq i$ and $1 \leq l,i \leq n$. Thus we can easily see that $r^{2}K=(e_{li})^{2}K=0$, but $rK=e_{li}K \neq 0$.
\end{itemize} 
Hence, in both cases $K \in$ $Nil_{M_n(R)}(M_{n}(M))$ and consequently, $Nil_{M_n(R)}(M_{n}(M)) = M_{n}(M)$\\

\begin{rmk}
	A Nil-semicommutative module need not necessarily be a Semicommutative module. The following is an example regarding the same.\\
\end{rmk}

\begin{exm}\label{eg3}
	Consider a left $R$-module $M$. Then, by Lemma \ref{lem1} $M_{n}(M)$ is Nil-semicommutative over $M_{n}(R)$ for $n \geq 2$. However, it is not  semicommutative  for $n \geq 4$ as shown below.\\
\end{exm}

\noindent\textbf{Proof:}~Let $A=\begin{pmatrix}
	0 & 1 & -1 & 0 &\dots & 0 \\ 
	0 & 0 & 0 & 0 &\dots & \vdots \\ 
	\vdots& \vdots & \vdots & \vdots & \ddots &  \vdots\\
	\vdots& \vdots & \vdots & \vdots & \ddots &  \vdots\\
	0 & 0 &  0 & 0 & \cdots & 0
\end{pmatrix}   \in M_{n}(R)$
and $K=\begin{pmatrix}
	0 & \cdots & 0 & 0 \\ 
	0 & \cdots & 0 & m \\ 
	0 & \cdots & 0 &  m\\
	0 & \cdots & 0 & 0\\
	\vdots & \vdots & \vdots & \vdots \\
	0 & \cdots & 0 & 0\\
\end{pmatrix} \in M_{n}(M)$ where $m \neq 0$.\\\\
Then, we can see that \\\\
$ AK = \begin{pmatrix}
	0 & 1 & -1 & 0 &\dots & 0 \\ 
	0 & 0 & 0 & 0 &\dots & \vdots \\ 
	\vdots& \vdots & \vdots & \vdots & \ddots &  \vdots\\
	\vdots& \vdots & \vdots & \vdots & \ddots &  \vdots\\
	0 & 0 &  0 & 0 & \cdots & 0
\end{pmatrix} \begin{pmatrix}
	0 & \cdots & 0 & 0 \\ 
	0 & \cdots & 0 & m \\ 
	0 & \cdots & 0 &  m\\
	0 & \cdots & 0 & 0\\
	\vdots & \vdots & \vdots & \vdots \\
	0 & \cdots & 0 & 0
\end{pmatrix} = 0 $ \\
\vspace{2 em}
\noindent But, for $ L = \begin{pmatrix}
	0 & 0 & 0 & 0 & \cdots & 0 \\ 
	0 & 0 & 1 & 0 & \cdots & 0 \\ 
	0 & 0 & 0 & 0 & \cdots & 0 \\
	\vdots & \vdots & \vdots & \vdots & \vdots & \vdots\\ 
	0 & 0 & 0 & 0 & 0 & 0
\end{pmatrix} \in M_n(R) $, we have \\\\
$ ALK = \begin{pmatrix}
	0 & 1 & -1 & 0 &\dots & 0 \\ 
	0 & 0 & 0 & 0 &\dots & \vdots \\ 
	\vdots& \vdots & \vdots & \vdots & \ddots &  \vdots\\
	\vdots& \vdots & \vdots & \vdots & \ddots &  \vdots\\
	0 & 0 &  0 & 0 & \cdots & 0
\end{pmatrix}  \begin{pmatrix}
	0 & 0 & 0 & 0 & \cdots & 0 \\ 
	0 & 0 & 1 & 0 & \cdots & 0 \\ 
	0 & 0 & 0 & 0 & \cdots & 0 \\
	\vdots & \vdots & \vdots & \vdots & \vdots & \vdots\\ 
	0 & 0 & 0 & 0 & 0 & 0
\end{pmatrix} \begin{pmatrix}
	0 & \cdots & 0 & 0 \\ 
	0 & \cdots & 0 & m \\ 
	0 & \cdots & 0 &  m\\
	0 & \cdots & 0 & 0\\
	\vdots & \vdots & \vdots & \vdots \\
	0 & \cdots & 0 & 0
\end{pmatrix} = \begin{pmatrix}
	0 & \cdots & 0 & m \\ 
	0 & \cdots & 0 & 0 \\ 
	0 & \cdots & \vdots & \vdots \\ 
	0 & 0 & 0 & 0
\end{pmatrix}  \neq 0 $\\\\
Hence, $_{M_n(R)}M_n(M)$ is not Semicommutative for $n\geq 4$.\\

%\begin{rmk}
%	For a field $F$, every left $F$- module $_FM$ is Semicommutative, Nilsemicommutative and Weakly Semicommutative. 
%\end{rmk}

\begin{prop}
	For a torsion-free module $_RM$, Nil$_R(M)$ = ${0}$.\\
\end{prop}

\noindent\textbf{Proof:} Let $m \in Nil_R(M)$.~Then, $\exists~ r_0 \in R$ and   $k_0 \in \mathbb{N}$ such that $r_0^{k_0}.m = 0 ~;~ r_0.m \neq 0$.~But, M is torsion-free so $m = 0$.~So, $ Nil_R(M) = 0$.\\

\begin{prop}
	Let $M$ be left $R$-module. If $M$ is torsion-free, then the following statements hold and are also equivalent.\\
	(i) $_RM$ is Semicommutative\\
	(ii) $_RM$ is Nil-semicommutative \\
	(iii) $_RM$ is Weakly semicommutative\\
\end{prop}

\noindent	\textbf{Proof:}
The proof follows easily from the fact that if $M$ is torsion free, then $Nil_{R}(M)=\{0\}$. \\ 

\noindent$(i) \Rightarrow (ii)$ Let $a \in R, m \in M$ be such that $am \in Nil_R(M)$. But, $M$ is torsion-free, so $Nil_R(M) = 0$ and so, $am = 0$. Thus, $aRm = 0 \subseteq Nil_R(M)$. Hence, $M$ is Nil-semicommutative module.\\ 

\noindent$(ii) \Rightarrow (iii)$
This follows from definition.\\ 

\noindent$(iii) \Rightarrow (i)$
Let $a \in R$  and $m \in M$ be such that $am = 0$.~Then, as $_RM$ is Weakly semicommutative module, $aRm \subseteq Nil_R(M) = 0$. So, $aRm = 0$. Hence, $_RM$ is Semicommutative module.\\

\begin{rmk}
	Nil modules are Semicommutative, Nil-semicommutative and also Weakly semicommutative.\\
\end{rmk}

%\begin{prop} Let $F$ be a field and $M$ be a left $F$-module.~Then, M is Semicommutative module $\Rightarrow$ $M$ is Nil semicommutative module.
%\end{prop}
%
%\textbf{proof} Let $a \in  F$, $m \in M$ satisfy $am \in Nil_F(M)$.~Now, since $am \in Nil_F(M)$ we have,~either $am = 0$ or $am \neq 0$.~If $am = 0$, then $aFm = 0$ $\subseteq  Nil_F(M)$.~If $am \neq 0$, then $\exists~ r \in F$  and $k \in N$ s.t. $r^k(am) = 0~;~ r(am) \neq 0$ $\Rightarrow$ $(r^ka).m = 0~;~ r(am) \neq 0$ $\Rightarrow$ $p.m = 0$, where $ p = r^ka$ $\Rightarrow$ $pFm = 0$ $\Rightarrow$ $r^katm = 0 ~ \forall~ t \in  F$.~Now, $r(atm) \neq 0$ as $r(atm) = 0 \Rightarrow (rat)m = 0 \Rightarrow (rat)Fm = 0 \Rightarrow (rat)t^{-1}m = 0 \Rightarrow ram = 0$ (which is absurd). So, $atm \in Nil_F(M) \Rightarrow aFm \subseteq Nil_F(M)$.~Hence, $M$ is Nil semicommutative module. \\

%	\begin{exm}
	
	%\end{exm}
	
	\begin{rmk}
		Submodules of Nil-semicommutative modules are Nil-semicommutative.\\
		%		but not so in case of direct sums and direct products.\\
	\end{rmk}
	
	%	The fact that the class of Nil-semicommutative modules is not closed under direct sums (or direct products) can be easily established from the following example in which we have shown that for two modules $M_1$ and $M_2$ over a ring $R$, $Nil_R(M_1 \oplus M_2) = Nil_R(M_1) \oplus Nil_R(M_2)$ may not hold in general. \\
	
	%	\begin{exm}
		%		
		%		Consider the $\mathbb{Z}$ modules $M_1 = ~\mathbb{Z}_2$ and $M_2 = ~\mathbb{Z}_4$. Then,   $ Nil_\mathbb{Z}(\mathbb{Z}_2) = \{ \bar{0}, \bar{1} \}$ and $ Nil_\mathbb{Z}(\mathbb{Z}_4)=\{ \bar{0}, \bar{1}, \bar{3} \}$. So,  $ Nil_\mathbb{Z}(\mathbb{Z}_2) \oplus  Nil_\mathbb{Z}(\mathbb{Z}_4)$ = $ \{ (\bar{0},\bar{0}), (\bar{0},\bar{1}), (\bar{0}, \bar{3}) \}$. However, we have $ \mathbb{Z}_2 \oplus ~ \mathbb{Z}_4 = \{ (\bar{0},\bar{0}), (\bar{0}, \bar{1}), (\bar{0}, \bar{2}), (\bar{0}, \bar{3}), (\bar{1}, \bar{0}), (\bar{1}, \bar{1}), (\bar{1}, \bar{2}), (\bar{1}, \bar{3}) \}$. So, $Nil_\mathbb{Z}(\mathbb{Z}_2 \times~ \mathbb{Z}_4) = \{ (\bar{0}, \bar{0}), (\bar{0}, \bar{1}), (\bar{0}, \bar{3}),
		%		(\bar{1}, \bar{0}), (\bar{1}, \bar{1}), (\bar{1}, \bar{3}) \}$.\\
		%		
		%		Hence, as we can see $ Nil_\mathbb{Z}(\mathbb{Z}_2) \oplus  Nil_\mathbb{Z}(\mathbb{Z}_4) \neq Nil_\mathbb{Z}(\mathbb{Z}_2 \oplus~ \mathbb{Z}_4)$. \\
		%	\end{exm}
	%	
	%	As a result, we can conclude that $M_1$ and $M_2$ being {\it Nil-semicommutative} may not imply that $M_1 \oplus M_2$ is also Nil-semicommutative.\\
	
	\begin{prop} For a ring $R$ and a left $R$ - module $M$, the following conditions are equivalent.\\
		(i) ~$M$ is Nil-semicommutative module.\\
		(ii) ~Every submodule of $_RM$ is Nil-semicommutative module.\\
		(iii) Every finitely generated submodule of $_RM$ is Nil semicommutative.\\
		(iv) Every cyclic submodule of $_RM$ is Nil-semicommutative module.\\
	\end{prop}
	
	\noindent\textbf{Proof:} $(i) \Rightarrow (ii)$ ~Let $N \subseteq M$ be a submodule and let $a$ $\in R$ , $n \in N$ be such that $an \in Nil_R(N)$.Then, $an \in Nil_R(M)$. So, $aRn \subseteq Nil_R(M)$. Thus, $\exists~ r_0 \in R$ and $k_0 \in \mathbb{N}$ such that $r_0^{k_0}(atn) = 0, r_0(atn) \neq 0 ~\forall~ t \in R$. But, $atn \in N$, so we have $atn \in Nil_R(N)$. Hence, $aRn \subseteq Nil_R(N)$. Thus, $N$ is Nil-semicommutative module.\\
	
	\noindent$(ii) \Rightarrow (iii)$ and $(iii) \Rightarrow (iv)$ ~follow clearly.\\  
	
	\noindent$(iv) \Rightarrow (i)$ Let $a \in R$ and $m \in M$  be such that $am \in Nil_R(M)$. This means that $\exists~ r_0 \in R$ and $k_0 \in \mathbb{N}$ such that $r_0^{k_0}(am) = 0~;~ r_0(am) \neq 0$. Now,~  $m=1.m \in Rm \leq M$. This implies that $ am \in Rm$ and consequently, $am \in Nil_R(Rm)$.~But, as $Rm$ is cyclic submodule of $M$, we have $aRm \subseteq Nil_R(Rm) \subseteq Nil_R(M)$. So, we have $aRm \subseteq Nil_R(M)$. Hence, $M$ is Nil-semicommutative.\\

	%		 Let $a \in R$ and $m \in M$  be such that $am \in Nil_R(M)$. This means that $\exists~ r_0 \in R$ and $k_0 \in \mathbb{N}$ such that $r_0^{k_0}(am) = 0~;~ r_0(am) \neq 0$. Now,~  $m=1.m \in Rm \leq M$. Now, $Rm$ being a submodule of $M$, implies that $ am \in Rm$ and consequently, $am \in Nil_R(Rm) \subseteq Nil_R(M)$. So, we have $aRm \subseteq Nil_R(M)$. Hence, $M$ is Nil semicommutative.\\ 
	
	Although we have seen in Lemma \ref{lem1} that for any given $_RM$ module and $n \geq 2$, $_{M_{n}(R)}M_n(M)$ is always a nilmodule and hence, a Nil-semicommutative  module. However, this may not be so in case of $_{T_{n}(R)}T_n(M)$. The following is an example in this regard.\\

	\begin{exm}\label{eg4}
		Consider the $\mathbb{Z}_{p^n}$ - module $\mathbb{Z}_{p^n}$. Then, for $n \geq 2$, $_{T_n(\mathbb{Z}_{p^n})}T_n(\mathbb{Z}_{p^n})$ is not Nil- semicommutative.\\
	\end{exm}
	
	\noindent	\textbf{Proof:} Let $A=\begin{pmatrix}
		\bar{1} & \bar{0} & \dots & \bar{0} \\ 
		\bar{0} & \bar{0} & \dots & \vdots \\ 
		\vdots& \vdots & \ddots &  \vdots\\
		\bar{0} & \bar{0} & \cdots & \bar{0}
	\end{pmatrix} $ and 	$ K = \begin{pmatrix}
		\bar{1} & \bar{0} & \dots & \bar{0} \\ 
		\bar{0} & \bar{0} & \dots & \vdots \\ 
		\vdots& \vdots & \ddots &  \vdots\\
		\bar{0} & \bar{0} & \cdots & \bar{0}
	\end{pmatrix} ~\in~ _{T_n(\mathbb{Z}_{p^n})}T_n(\mathbb{Z}_{p^n}) $\\\\
	
	\noindent Then, $AK = \begin{pmatrix}
		\bar{1} & \bar{0} & \dots & \bar{0} \\ 
		\bar{0} & \bar{0} & \dots & \vdots \\ 
		\vdots& \vdots & \ddots &  \vdots\\
		\bar{0} & \bar{0} & \cdots & \bar{0}
	\end{pmatrix}  \in Nil_{T_n(\mathbb{Z}_{p^n})}(T_n(\mathbb{Z}_{p^n}))$. \\\\
	
	\noindent as $ \exists~ P = \begin{pmatrix}
		\bar{p} & \bar{0} & \dots & \bar{0} \\ 
		\bar{0} & \bar{0} & \dots & \vdots \\ 
		\vdots& \vdots & \ddots &  \vdots\\
		\bar{0} & \bar{0} & \cdots & \bar{0}
	\end{pmatrix} ~\in~ T_n(\mathbb{Z}_{p^n}) $ such that 
	$P^n(AK) = 0$ but $ P(AK) \neq 0$ \\
	
	\noindent But, for $ L = \begin{pmatrix}
		\bar{p}^{n-1} & \bar{0} & \dots & \bar{0} \\ 
		\bar{0} & \bar{0} & \dots & \vdots \\ 
		\vdots& \vdots & \ddots &  \vdots\\
		\bar{0} & \bar{0} & \cdots & \bar{0}
	\end{pmatrix}$, we claim that $ ALK = \begin{pmatrix}
		\bar{p}^{n-1} & \bar{0} & \dots & \bar{0} \\ 
		\bar{0} & \bar{0} & \dots & \vdots \\ 
		\vdots& \vdots & \ddots &  \vdots\\
		\bar{0} & \bar{0} & \cdots & \bar{0}
	\end{pmatrix} \\
	\vspace{0.5 cm} 
	\notin Nil_{T_n(\mathbb{Z}_{p^n})}(T_n(\mathbb{Z}_{p^n}))$.\\
	
	\noindent For this, suppose $\exists~ S = \begin{pmatrix}
		\bar{s_1} & * & \dots & * \\ 
		\bar{0} & \bar{s_2} & \dots & * \\ 
		\vdots& \vdots & \ddots &  \vdots\\
		\bar{0} & \bar{0} & \cdots & \bar{s_n}
	\end{pmatrix} ~\in~ T_n(\mathbb{Z}_{p^n}) $ such that 
	$S^k(ALK) = 0$ \\ but $ S(ALK) \neq 0$. \\
	
	\noindent	Now, $S^k(ALK) = 0$\\\\
	$\implies \begin{pmatrix}
		\bar{s_1} & * & \dots & * \\ 
		\bar{0} & \bar{s_2} & \dots & * \\ 
		\vdots& \vdots & \ddots &  \vdots\\
		\bar{0} & \bar{0} & \cdots & \bar{s_n}
	\end{pmatrix}^k \begin{pmatrix}
		\bar{p}^{n-1} & \bar{0} & \dots & \bar{0} \\ 
		\bar{0} & \bar{0} & \dots & \vdots \\ 
		\vdots& \vdots & \ddots &  \vdots\\
		\bar{0} & \bar{0} & \cdots & \bar{0}
	\end{pmatrix} = 0$\\\\ 
	
	\noindent $\implies \begin{pmatrix}
		\bar{s_1}^k\bar{p}^{n-1} & * & \dots & * \\ 
		\bar{0} & \bar{0} & \dots & * \\ 
		\vdots& \vdots & \ddots &  \vdots\\
		\bar{0} & \bar{0} & \cdots & \bar{0}
	\end{pmatrix} = 0 $\\\\
	
	\noindent	 $\implies {s_1}^kp^{n-1} = p^ns$ for some $s \in \mathbb{Z}$.\\
	
	\noindent	$\implies p~|~ s_1^k $\\
	
	\noindent Consequently, $ p~|~s_1$ and thus, $s_1 = pl$ for some $ l \in \mathbb{Z}$.\\
	
	\noindent So,  ~$ S(ALK) = ~\begin{pmatrix}
		\bar{pl} & * & \dots & * \\ 
		\bar{0} & \bar{a_2} & \dots & * \\ 
		\vdots& \vdots & \ddots &  \vdots\\
		\bar{0} & \bar{0} & \cdots & \bar{a_n}
	\end{pmatrix}\begin{pmatrix}
		\bar{p}^{n-1} & \bar{0} & \dots & \bar{0} \\ 
		\bar{0} & \bar{0} & \dots & \vdots \\ 
		\vdots& \vdots & \ddots &  \vdots\\
		\bar{0} & \bar{0} & \cdots & \bar{0}
	\end{pmatrix}  ~=~ 0$, which is a contradiction.\\
	
	\noindent Therefore, there does not exist such $S$ and so, $ ALK \notin Nil_{T_n(\mathbb{Z}_{p^n})}(T_n(\mathbb{Z}_{p^n}))$. Hence, $_{T_n(\mathbb{Z}_{p^n})}T_n(\mathbb{Z}_{p^n})$ is not Nilsemicommutative.\\\\
	
	Moreover, unlike in Weakly semicommutative case as in \textcolor{blue}{\cite{1}}, $_RM$ being a Nil-semicommutative module may not necessarily imply that $_{T_n(R)}T_n(M)$ is also Nil-semicommutative. The following is an example in this regard.\\\\
	
	\begin{exm}\label{eg5}
		Consider the $\mathbb{Z}_p$- module $\mathbb{Z}_p$. Then, $\mathbb{Z}_p$ is Nil-semicommutative as it is torsion free. However, for $ n\geq 2$, $_{T_n(\mathbb{Z}_p)}T_n(\mathbb{Z}_p)$ is not Nil-semicommutative.\\
	\end{exm}
	
	\noindent \textbf{Proof:} Let $A = \begin{pmatrix}
		\bar{p-1} & \bar{0} & \dots & \bar{0} \\ 
		\bar{0} & \bar{p-1} & \dots & \vdots \\ 
		\vdots& \vdots & \ddots &  \vdots\\
		\bar{0} & \bar{0} & \cdots & \bar{p-1}
	\end{pmatrix} $ and 	$ K = \begin{pmatrix}
		\bar{1} & \bar{0} & \dots & \bar{0} \\ 
		\bar{0} & \bar{1} & \dots & \vdots \\ 
		\vdots& \vdots & \ddots &  \vdots\\
		\bar{0} & \bar{0} & \cdots & \bar{1}
	\end{pmatrix} ~\in~ _{T_n(\mathbb{Z}_p)}T_n(\mathbb{Z}_p) $\\\\
	
	\noindent Then, $ AK = \begin{pmatrix}
		\bar{p-1} & \bar{0} & \dots & \bar{0} \\ 
		\bar{0} & \bar{p-1} & \dots & \vdots \\ 
		\vdots& \vdots & \ddots &  \vdots\\
		\bar{0} & \bar{0} & \cdots & \bar{p-1}
	\end{pmatrix} $ $ \in  Nil_{T_n(\mathbb{Z}_p)}(T_n(\mathbb{Z}_p)) $\\\\
	
	\noindent as $ \exists~ P = \begin{pmatrix}
		\bar{0} & \bar{0} & \dots & \bar{1} \\ 
		\bar{0} & \bar{0} & \dots & \vdots \\ 
		\vdots& \vdots & \ddots &  \vdots\\
		\bar{0} & \bar{0} & \cdots & \bar{0}
	\end{pmatrix} ~\in~ _{T_n(\mathbb{Z}_p)}T_n(\mathbb{Z}_p) $ such that 
	$P^2(AK) = 0$ but $P(AK) \neq 0$.\\\\
	
	\noindent But, for $ L = \begin{pmatrix}
		\bar{0} & \bar{0} & \dots & \bar{p-1} \\ 
		\bar{0} & \bar{0} & \dots & \vdots \\ 
		\vdots& \vdots & \ddots &  \vdots\\
		\bar{0} & \bar{0} & \cdots & \bar{0}
	\end{pmatrix} ~\in~ T_n(\mathbb{Z}_p) $, we claim that \\\\
	
	\noindent $ ALK = \begin{pmatrix}
		\bar{p-1} & \bar{0} & \dots & \bar{0} \\ 
		\bar{0} & \bar{p-1} & \dots & \vdots \\ 
		\vdots& \vdots & \ddots &  \vdots\\
		\bar{0} & \bar{0} & \cdots & \bar{p-1}
	\end{pmatrix} \begin{pmatrix}
		\bar{0} & \bar{0} & \dots & \bar{p-1} \\ 
		\bar{0} & \bar{0} & \dots & \vdots \\ 
		\vdots& \vdots & \ddots &  \vdots\\
		\bar{0} & \bar{0} & \cdots & \bar{0}
	\end{pmatrix}
	\begin{pmatrix}
		\bar{1} & \bar{0} & \dots & \bar{0} \\ 
		\bar{0} & \bar{1} & \dots & \vdots \\ 
		\vdots& \vdots & \ddots &  \vdots\\
		\bar{0} & \bar{0} & \cdots & \bar{1}
	\end{pmatrix} $
	\vspace{2mm}
	
	\noindent $ = \begin{pmatrix}
		\bar{0} & \bar{0} & \dots & (\bar{p-1})^2 \\ 
		\bar{0} & \bar{0} & \dots & \vdots \\ 
		\vdots& \vdots & \ddots &  \vdots\\
		\bar{0} & \bar{0} & \cdots & \bar{0}
	\end{pmatrix} $
	\vspace{2mm}
	
	\noindent $ = \begin{pmatrix}
		\bar{0} & \bar{0} & \dots & \bar{1} \\ 
		\bar{0} & \bar{0} & \dots & \vdots \\ 
		\vdots& \vdots & \ddots &  \vdots\\
		\bar{0} & \bar{0} & \cdots & \bar{0}
	\end{pmatrix} \notin Nil_{T_n(\mathbb{Z}_p)}(T_n(\mathbb{Z}_p)) $\\
	
	\noindent For this, suppose $\exists~ S = \begin{pmatrix}
		\bar{s_1} & * & \dots & * \\ 
		\bar{0} & \bar{s_2} & \dots & * \\ 
		\vdots& \vdots & \ddots &  \vdots\\
		\bar{0} & \bar{0} & \cdots & \bar{s_n}
	\end{pmatrix} ~\in~ T_n(\mathbb{Z}_{p}) $ such that 
	$S^k(ALK) = 0$ but $ S(ALK) \neq 0$.\\
	
	\noindent Then, as in the previous Example \ref{eg4} we can easily prove that $ s_1 = pl$ for some $ l \in \mathbb{Z}$. This will consequently give us $S(ALK) = 0$, thereby contradicting our assumption.Therefore, there does not exist such $S$ and so, $ ALK \notin Nil_{T_n(\mathbb{Z}_p)}(T_n(\mathbb{Z}_p)) $. Hence, $_{T_n(\mathbb{Z}_p)}T_n(\mathbb{Z}_p)$ is not Nil-semicommutative.\\
	
	Further, recall from \textcolor{blue}{\cite{lee}} that for a given left $R$-module $M$ and $ A = (a_{ij}) \in M_n(R)$, the ring of  $n \times n$ matrices over $R$, $AM = \{ (a_{ij}m): m \in M \}$.~Moreover, for $ V = \sum_{i=1}^{n-1} E_{i,i+1}, n \geq 2$ where $E_{i,j}$ are matrix units, $V_n(R) = RI_n + RV + RV^2 + \hdots + RV^{n-1}$ forms a ring and $ V_n(M) = MI_n + MV + MV^2 + \hdots + MV^{n-1}$ forms a left module over the ring $V_n(R)$ . In addition to this, there exists a ring isomorphism $\theta : V_n(R) \rightarrow \frac{R[x]}{(x^n)} $ ~defined by~ $ \theta (r_0I_n + r_1V + r_2V^2 + \hdots + r_{n-1}V^{n-1}) = r_0 + r_1x + \hdots r_{n-1}x^{n-1} + (x^n)$ and an abelian group isomorphism $\phi : V_n(M) \rightarrow \frac{M[x]}{(M[x](x^n))}$ ~defined by~ $\phi(m_0 + m_1V + m_2V^2 + \hdots + m_1V^{n-1}) = m_0 + m_1x + m_2x^2 + \hdots + m_{n-1}x^{n-1} + M[x](x^n)$ such that~ $ \phi (AW) = \theta(A) \phi(W) ~\forall~ A \in V_n(R), W \in V_n(M)$. \\
	
	\begin{rmk}
		For a left $R$ - module $M$, $_RM$ is Nil-semicommutative may not imply that $_{V_n(R)}V_n(M)$ is Nil-semicommutative. We can consider the following  Example \ref{eg6} in this regard.\\
	\end{rmk}
	
	\begin{exm}\label{eg6}
		Consider the $\mathbb{Z}_p$- module $\mathbb{Z}_p$. Then, $_{\mathbb{Z}_p}\mathbb{Z}_p$ is Nil-semicommutative. However, for $n\geq 2$, $_{V_n(\mathbb{Z}_p)}V_n(\mathbb{Z}_p)$ is not Nil-semicommutative.\\
	\end{exm}
	
	\noindent \textbf{Proof:} Let $A = \begin{pmatrix}
		\bar{p-1} & \bar{0} & \dots & \bar{0} \\ 
		\bar{0} & \bar{p-1} & \dots & \vdots \\ 
		\vdots& \vdots & \ddots &  \vdots\\
		\bar{0} & \bar{0} & \cdots & \bar{p-1}
	\end{pmatrix} $ and 	$ K = \begin{pmatrix}
		\bar{1} & \bar{0} & \dots & \bar{0} \\ 
		\bar{0} & \bar{1} & \dots & \vdots \\ 
		\vdots& \vdots & \ddots &  \vdots\\
		\bar{0} & \bar{0} & \cdots & \bar{1}
	\end{pmatrix} ~\in~ _{V_n(\mathbb{Z}_p)}V_n(\mathbb{Z}_p) $\\\\
	
	\noindent Then, similarly as in Example \ref{eg5}, $ AK = \begin{pmatrix}
		\bar{p-1} & \bar{0} & \dots & \bar{0} \\ 
		\bar{0} & \bar{p-1} & \dots & \vdots \\ 
		\vdots& \vdots & \ddots &  \vdots\\
		\bar{0} & \bar{0} & \cdots & \bar{p-1}
	\end{pmatrix} $ $ \in  Nil_{V_n(\mathbb{Z}_p)}(V_n(\mathbb{Z}_p)) $\\\\
	
	%    \noindent as $ \exists~ P = \begin{pmatrix}
		%    	\bar{0} & \bar{0} & \dots & \bar{1} \\ 
		%    	\bar{0} & \bar{0} & \dots & \vdots \\ 
		%    	\vdots& \vdots & \ddots &  \vdots\\
		%    	\bar{0} & \bar{0} & \cdots & \bar{0}
		%    \end{pmatrix} ~\in~ _{V_n(\mathbb{Z}_p)}V_n(\mathbb{Z}_p) $ such that 
	%    $P^2(AK) = 0$ but $P(AK) \neq 0$.\\\\
	
	\noindent But, for $ L = \begin{pmatrix}
		\bar{0} & \bar{0} & \dots & \bar{p-1} \\ 
		\bar{0} & \bar{0} & \dots & \vdots \\ 
		\vdots& \vdots & \ddots &  \vdots\\
		\bar{0} & \bar{0} & \cdots & \bar{0}
	\end{pmatrix} ~\in~ V_n(\mathbb{Z}_p) $, claim \\\\
	
	\noindent $ ALK = \begin{pmatrix}
		\bar{p-1} & \bar{0} & \dots & \bar{0} \\ 
		\bar{0} & \bar{p-1} & \dots & \vdots \\ 
		\vdots& \vdots & \ddots &  \vdots\\
		\bar{0} & \bar{0} & \cdots & \bar{p-1}
	\end{pmatrix} \begin{pmatrix}
		\bar{0} & \bar{0} & \dots & \bar{p-1} \\ 
		\bar{0} & \bar{0} & \dots & \vdots \\ 
		\vdots& \vdots & \ddots &  \vdots\\
		\bar{0} & \bar{0} & \cdots & \bar{0}
	\end{pmatrix}
	\begin{pmatrix}
		\bar{1} & \bar{0} & \dots & \bar{0} \\ 
		\bar{0} & \bar{1} & \dots & \vdots \\ 
		\vdots& \vdots & \ddots &  \vdots\\
		\bar{0} & \bar{0} & \cdots & \bar{1}
	\end{pmatrix} $
	\vspace{2mm}
	
	%    \noindent $ = \begin{pmatrix}
		%    	\bar{0} & \bar{0} & \dots & (\bar{p-1})^2 \\ 
		%    	\bar{0} & \bar{0} & \dots & \vdots \\ 
		%    	\vdots& \vdots & \ddots &  \vdots\\
		%    	\bar{0} & \bar{0} & \cdots & \bar{0}
		%    \end{pmatrix} $
	%    \vspace{2mm}
	
	\noindent $ = \begin{pmatrix}
		\bar{0} & \bar{0} & \dots & \bar{1} \\ 
		\bar{0} & \bar{0} & \dots & \vdots \\ 
		\vdots& \vdots & \ddots &  \vdots\\
		\bar{0} & \bar{0} & \cdots & \bar{0}
	\end{pmatrix} \notin Nil_{V_n(\mathbb{Z}_p)}(V_n(\mathbb{Z}_p)) $\\
	
	\noindent For this, suppose $\exists~ S = \begin{pmatrix}
		\bar{s} & * & \dots & * \\ 
		\bar{0} & \bar{s} & \dots & * \\ 
		\vdots& \vdots & \ddots &  \vdots\\
		\bar{0} & \bar{0} & \cdots & \bar{s}
	\end{pmatrix} ~\in~ V_n(\mathbb{Z}_{p}) $ such that 
	$S^k(ALK) = 0$ but $ S(ALK) \neq 0$.\\
	
	\noindent Then, similarly as in the Example \ref{eg4} we can prove that $ s = pl$ for some $ l \in \mathbb{Z}$. Consequently, $S(ALK) = 0$, which is a contradiction.Therefore, $\nexists$ such $S$ and so, $ ALK \notin Nil_{V_n(\mathbb{Z}_p)}(V_n(\mathbb{Z}_p)) $. Hence, $_{V_n(\mathbb{Z}_p)}V_n(\mathbb{Z}_p)$ is not Nil-semicommutative.\\
	\vspace{2mm}

	In \textcolor{blue}{\cite{1}} (Proposition $3.15$), it was proved that if $_RM$ is Weakly semicommutative, then for $n\geq 2$ $ \frac{M[x]}{(M[x](x^n))}$ is also a Weakly semicommutative module over $\frac{R[x]}{(x^n)} $. However, as stated in the following remark, this may not be true in the Nil-semicommutative case.\\
	
	\begin{rmk}
		For a left $R$- module $M$, $_RM$ is Nil-semicommutative may not imply that $ \frac{M[x]}{(M[x](x^n))}$ is  Nil-semicommutative over $\frac{R[x]}{(x^n)} $. For this, as we have seen in Example \ref{eg6}, it does not hold true in  the $_{V_n(R)}V_n(M)$ case. Hence, it will not hold true in this case.\\
		
		%    	\textbf{Proof:} The proof follows from the fact that $\exists$ module isomorphism from $V_n(M)$ to $ \frac{M[x]}{(M[x](x^n))}$ and ring isomorphism from $V_n(R)$ to $\frac{R[x]}{(x^n)} $.
	\end{rmk}

	%Further, we prove the module counterparts for Proposition 1.1 and Proposition 1.2. However while doing so, we have also deduced that unlike the case in their ring counterparts, an SC module $M$ may not necessarily be an NSC module. While an example in support of the above has been discussed later in the paper, the conditions required for the generalisation to hold are discussed as follows.  

	We recall from \textcolor{blue}{\cite{ssev}} that for $S_0$ denoting the set of nonzero divisors of $R$, $T(M) = [ m \in M : rm = 0$ for some $r \in S_0]$. Also, for a module $M$ over a general ring $R$, $T(M) \subseteq Tor(M)$ and for $R$ domain, $T(M) = Tor(M)$. If $R$ is not domain, then the equality may not hold. For instance, in the ${\mathbb{Z}_6}$ - module  $\mathbb{Z}_6$, $\bar{2}$.$\bar{3}$ = $\bar{0}$ $\Rightarrow$ $\bar{3}$ $\in Tor(M)$ but $\bar{3}$ $\notin$ $T(M)$ as $\bar{3}$ is a zero divisor in $\mathbb{Z}_6$. Moreover, if $R$ is a commutaive domain then, $Tor(M)$ is a submodule of $M$. Further, Rege and Bhuphang in \textcolor{blue}{\cite{br}} proved that for a left $R$- module $M$, $_RM$ is Semicommutative module $\Rightarrow$ $T(M)$ is a submodule of $M$. In addition to this, Ansari and Herachandra in \textcolor{blue}{\cite{ssev}} generalised this result when they proved that for a Weakly semicommutative module $M$ over a domain $D$, $T(M)$ is a submodule of $M$. Similarly, we have continued the trend and proved that the result holds true in case of Nil-semicommutative modules as well.\\ 
	
	\begin{prop}
		For a Nil-semicommutative module $M$ over a domain D, $T(M)$ is a submodule of $M$.
	\end{prop}
	
	\noindent\textbf{Proof:} The proof follows easily from the fact that every Nil-semicommutative module is Weakly semicommutative.\\

	%	\begin{prop}
		%		For a torsion-free module M, we have
		%		$M$ is Nil-semicommutative module $\iff$ $M/Tor(M)$ is Nil semicommutative module.\\
		%	\end{prop} 
	
	\begin{rmk}
		For a torsion-free module $M$, the following are equivalent\\
		(i) $M$ is Nil-semicommutative module\\
		(ii) $M/Tor(M)$ is Nil-semicommutative module. \\

	\end{rmk}
	
	%	\begin{rmk}
		%		If $R$ is an Integral domain and $M$ is torsion-free module, then we have $M$ is Nil semicommutative module $\Rightarrow$  
		%	\end{rmk}
	%\begin{rmk}
	%For a given left $R$ module $_RM$, $Tor(M)$ is never Torsion-free.
	%\end{rmk}

	\begin{prop} Let $R$ be a commutative ring in which the degree of nilpotency for any $r \in Nil(R) ~ is > 2$.~Then, $R$ is Nil-semicommutative ring $\Rightarrow$ $_RR$ is Nil-semicommuatative module.\\
	\end{prop}
	
	\noindent\textbf{Proof:} ~Let a $\in R$ and $m \in R$ satisfy $am \in Nil_R(R)$. Then, two cases arise. For $am = 0$, the proof is easy.~For $am \neq 0$, $\exists~ r \in R$ and $k \in \mathbb{N}$ s.t. $r^k (am) = 0$ but $r(am) \neq 0$. So, for $p = r^ka$ we have, $ pm = 0$. This implies that $ pm \in Nil(R)$ and hence, $pRm \subseteq Nil(R)$. Consequently, we have $(r^ka)(t)m \in Nil(R)$.~Now, if $r^k(atm) = 0$, then we have $atm \in Nil_R(R)$ and so, $aRm \subseteq Nil_R(R)$.
	And if $r^k(atm) \neq 0$, then $\exists ~k_2 \in \mathbb{N}$~ such that $~{(r^k(atm))}^{k_2} = 0 $. Consequently, we have $r^{k.k_2}(atm)^{k_2}(atm)^{k_2(k-1)+1}=0$ and thus, $(ratm)^{kk_2}(atm) = 0 $. Now, for $r_* = ratm, q = kk_2 \in \mathbb{N}$, we have $(r_*)^q (atm) = 0~$ .We then claim that $r_*(atm) \neq 0.$ For if $r_*(atm) = 0$ then, $(ratm)(atm) = r(atm)^2 = 0$ which implies that $ r^{2k-1}.r(atm)^2 = 0$ and so, $ (r^k(atm))^2 = 0$, which is a contradiction.Therefore, in both  cases, $atm \in Nil_R(R)$ $\Rightarrow aRm \subseteq Nil_R(R) $. Hence,~$_RR$ is Nil-semicommutative module.\\

	\begin{prop}
		Let $\phi: R \rightarrow R'$ be a ring homomorphism.~Then for a left $R'$- module $M$, $M$ is also a left $R$ - module via $rm = \phi(r)m$.~Moreover, if $\phi$ is onto then, the following are equivalent.\\
		(i) $M$ is Nil-semicommutative $R'$-module. \\
		(ii) $M$ is Nil-semicommutative $R$-module.\\
	\end{prop}
	
	\noindent\textbf{Proof:}
	Let $a$ $\in$ $R$, $m \in M$ such that $am \in Nil_R(M)$. Then, $\exists~ r \in R, n \in \mathbb{N}$ such that $r^k(am) = 0; r(am) \neq 0$. So, $(r^ka)m = 0~;~ (ra)m \neq 0$ and consequently, $\phi(r^ka)m = 0~;~ \phi(ra)m \neq 0$.This eventually implies that $\phi(r)^k\phi(a)m = 0~;~ \phi(r)\phi(a)m \neq 0 $ and hence, $\phi(a)m \in Nil_{R'}(M)$. But, $M$ is Nil-semicommutative $R'$ - module.~So, $\phi(a)R'm$ $\subseteq$ $Nil_{R'}(M)$ and  $\phi(a)r'_1m \in Nil_{R'}(M)~ \forall r'_1 \in R'$. But $\phi$ is onto so $\exists ~r_1 \in R$ such that $\phi(r_1) = r'_1$. So,  $\phi(a)\phi(r_1)m \in Nil_{R'}(M)$ and consequently, $\exists~ \phi(r_0) \in R', k \in \mathbb{N} $~s.t $\phi(r_0)^{k_0}(\phi(a)\phi(r_1)m) = 0~;~ \phi(r_0)(\phi(a)\phi(r_1)m) \neq 0$ Eventually, we then have $\phi(r_0^{k_0}ar_1)m = 0~;~ \phi(r_0ar_1)m \neq 0$ and hence, ${r_0}^{k_0}(ar_1m) = 0~;~ r_0(ar_1m) \neq 0$. Thus, we have $ ar_1m \in Nil_R(M)$ and consequently, $aRm \subseteq Nil_R(M)$.~Hence, $M$ is Nil-semicommutative $R$- module.\\
	
	Conversely, let $b \in R'$, $m \in M$ such that $bm \in Nil_{R'}(M)$.Then, $\exists ~r' \in R', k \in \mathbb{N}$ such that $(r')^kbm = 0 ~;~ r'(bm)$ $\neq$ 0.~Now, $\phi$ is onto so~ $\exists$ $a \in R$ and $r \in R$ such that $b = \phi(a)$ and $r'= \phi(r)$.~So, we have $\phi(r)^k\phi(a)m = 0 ~;~ \phi(r)(\phi(a)m) \neq 0$ and eventually, $\phi(r^ka)m = 0 ~;~ \phi(ra)m) \neq 0$.~This implies that $(r^ka)m = 0 ; (ra)m) \neq 0$ and consequently, $am \in Nil_R(M)$ and $aRm \subseteq Nil_R(M)$. Further, we also have $(r_0^{k_0}ar)m = 0~;~ (r_0ar)m \neq 0$ for some $r_0 \in R, k_0 \in \mathbb{N}$. This implies that $\phi(r_0^{k_0}ar)m = 0~;~ \phi(r_0ar)m \neq 0$ and eventually,
	$\phi(r_0)^{k_0}(\phi(a)\phi(r)m) = 0 ~;~ \phi(r_0)(\phi(a)\phi(r))m) \neq 0$. Thus, we have $br'm \in Nil_R'(M)$ and hence,~$M$ is Nil-semicommutative $R'$- module.\\

	\begin{prop}\label{prop 2.9}
		For a ring $R$ and a multiplicatively closed subset $S$ of $C(R) \backslash \{0\} $ containing $1$, $_RM$ is Nil-semicommutative module $\iff$ $S^{-1}M$ is Nil-semicommutative $S^{-1}R$ module.\\
	\end{prop}
	
	\noindent\textbf{Proof:} Let $\dfrac{r}{s} \in S^{-1}R, \dfrac{m}{q} \in S^{-1}M$ be such that $\dfrac{r}{s}.\dfrac{m}{q} \in Nil_{S^{-1}R}(S^{-1}M) $. Then, we have $(\dfrac{r_0}{s_0})^{k_0}(\dfrac{r}{s}.\frac{m}{q}) = 0 ~;~  (\dfrac{r_0}{s_0})(\dfrac{r}{s}.\dfrac{m}{q}) \neq 0 $ for some $ r_0 \in R, k_0 \in \mathbb{N}$. This implies that $u_1.r_0^{k_0}rm = 0$ for some $u \in S ~;~ u.r_0rm \neq 0 ~\forall~ u \in S$. Now, since $u_1 \in S \subseteq C(R)$, we have $u_1rm \in  Nil_R(M) $. Consequently, we have $u_1rlm \in Nil_R(M) ~\forall~ l \in R$. Therefore, $ \exists~ t_1 \in R, k_1 \in \mathbb{N}$ such that $t_1^{k_1}(u_1rlm) = 0 ~;~ t_1(u_1rlm) \neq 0$. This implies that for $s_1, s, p$ and $q \in S$, we have $u_1(t_1^{k_1}rlm.1 - 0.s_1spq) = 0 ~;~ u_1(t_1rlm.1 - 0.s_1spq) \neq 0$. Thus, we have  $\dfrac{t_1^{k_1}rlm}{s_1^{k1}spq} = 0$ and with an easy 'proof by contradiction' we also see that $\dfrac{t_1rlm}{s_1spq} \neq 0$. Hence, we have $\dfrac{r}{s}\dfrac{l}{p}\dfrac{m}{q} \in Nil_{S^{-1}R}(S^{-1}M)$ for arbitrary $\dfrac{l}{p} \in S^{-1}R$. And consequently, we have $_{S^{-1}R}S^{-1}(M)$ is Nil-semicommutative module.\\ 
	
	For the converse, let $a \in R, m \in M$ such that $am \in Nil_R(M)$. Then $ \exists ~ r_0 \in R, k_0 \in \mathbb{N}$ such that  $r_0^{k_0}am =0 ~;~ ram \neq 0$. This implies that for $s_0, s$ and $q \in S$, we have $u_1(r_0^{k_0}am.1 - 0.s_0^{k_0}sq) = 0$ for some $u_1 \in S$ ~;~ $u(ram.1 - 0.s_0sq) \neq 0 ~\forall~ u \in S$. Consequently, we have $(\dfrac{r_0}{s_0})^{k_0}(\dfrac{a}{s}.\dfrac{m}{q}) = 0 $~;~$(\dfrac{r_0}{s_0})(\dfrac{a}{s}.\dfrac{m}{q}) \neq 0$. Therefore, $\dfrac{a}{s}.\dfrac{m}{q} \in Nil_{S^{-1}R}(S^{-1}M)$ and consequently, $\dfrac{a}{s}(S^{-1}R)\dfrac{m}{q} \subseteq Nil_{S^{-1}R}(S^{-1}M)$. This implies that $\dfrac{a}{s}\dfrac{l}{p}\dfrac{m}{q} \in Nil_{S^{-1}R}(S^{-1}M) ~\forall~ \dfrac{l}{p} \in S^{-1}R $ and hence, we have $\dfrac{t_1^{k_1}alm}{s_1^{k1}spq} = 0$ and $\dfrac{t_1alm}{s_1^{k_1}spq} \neq 0$. Eventually, we then arrive at $u_1t_1alm = 0$ ~for some $u_1 \in S$ and $ut_1alm \neq 0 ~\forall~ u \in S$. Then, for $ u = 1$ and with the help of an easy proof, we finally have $ t_1^{k_1}alm = 0 ~;~ t_1alm \neq 0$. Therefore, $aRm \subseteq Nil_R(M)$ and thus, $_RM$ is Nil-semicommutative module. \\
	\vspace{2mm}

	We recall from \textcolor{blue}{\cite{lee}} that for a left $R$-module $M$, $ R[x,x^{-1}]$ is the ring of Laurent polynomials in $x$ with coefficients in $R$ and $ M[x,x^{-1}]$ forms a left module over the ring $ R[x,x^{-1}]$. Now, as a corollary to the above proposition, we deduce a relation between the $_{R[x]}M[x]$ module  and the $_{R[x,x^{-1}]}M[x,x^{-1}]$ module. \\
	
	\begin{corollary} \label{cor}
		For a left $R$-module $M$, $_{R[x]}M[x]$ is Nil-semicommutative $\iff _{R[x,x^{-1}]}M[x,x^{-1}]$ is Nil-semicommutative. \\
	\end{corollary}
	
	\noindent\textbf{Proof:} Let $ S = \{ 1,x,x^2,x^3,...\}$. Then $S$ is a multiplicatively closed subset of $C(R)\backslash \{0\}$ containing $1$.
	Now, for $ a_{-m}x^{-m} + a_{-(m-1)}x^{-(m-1)} + \cdots + a_{-1}x^{-1} + a_0 + a_1x + \cdots a_{-n}x^{n} \in R[x,x^{-1}]$ ~it is easy to see that
	$\dfrac{1}{x^m} [ {a_{-m} + a_{-(m-1)}x + \cdots + a_{-1}x^{m-1} + a_0x^m + a_1x^{m+1} + \cdots a_{-n}x^{m+n}}] \in S^{-1}R[x]$.~So, $R[x,x^{-1}] \subseteq S^{-1}R[x]$. Similarly, $ S^{-1}R[x] \subseteq R[x,x^{-1}]$.~Therefore, $S^{-1}R[x] = R[x,x^{-1}]$ and in the same way, $ S^{-1}M[x] = M[x,x^{-1}]$. Hence, by Proposition \ref{prop 2.9} the result follows clearly.\\
	
	\section{Concluding remark}
	We conclude this note with the following questions.
	\begin{quest}
		For a Nil-semicommutative module, $M$ and its submodule $N$, does it imply that $M/N$ is Nil-semicommutative? If so, what is the required condition?
	\end{quest} 
	\begin{quest}
		Under what condition will $_RR$ being a Nil-semicommutative module imply that $R$ is Nil-semicommutative?
	\end{quest} 
	\begin{quest}
		What is the relation between $Nil_R(M[x])$ and $Nil_R(M)[x]$?\\\\
	\end{quest}

\bibliography{sn-bibliography}% common bib file

\end{document}